\documentclass[12pt]{amsart}
\usepackage[top=1in, bottom=1in, left=1in, right=1in]{geometry}
\usepackage{amsfonts}
\usepackage{amsmath}
\usepackage{comment}
\usepackage{amssymb}
\usepackage{graphicx}
\usepackage{bbm}
\usepackage{comment}
\usepackage{mathrsfs}
\numberwithin{equation}{section}
\usepackage{times}
\usepackage[backend=biber, sorting=nyt, url=false,
maxnames = 100,
doi=false]{biblatex}
\usepackage{siunitx}
\usepackage[usenames,dvipsnames]{color}
\usepackage{comment}
\usepackage{xcolor}
\usepackage{mathtools}
\usepackage{bm}
\usepackage{esvect}
\usepackage{hyperref}

\hypersetup{
    colorlinks = true,
linkcolor={black},
urlcolor={blue},
citecolor={blue},    
urlcolor = {blue},
citebordercolor = {0.33 .58 0.33},
 linkbordercolor = {0.99 .28 0.23},
 breaklinks=true}
 
\newcommand{\R}{\mathbb{R}}
\newcommand{\Q}{\mathbb{Q}}
\newcommand{\N}{\mathbb{N}}
\newcommand{\Z}{\mathbb{Z}}

\newtheorem{thm}{Theorem}[section]

\newtheorem{rem}[thm]{Remark}
\newtheorem{conj}[thm]{Conjecture}

\theoremstyle{remark}

\usepackage{geometry}
\title{Erd\H{o}s-Moser Equation in Arithmetic Progressions}
\author{Anji Dong, Vi Anh Nguyen, Alexandru Zaharescu}

\address{
Anji Dong: Department of Mathematics,
University of Illinois Urbana-Champaign,
Altgeld Hall, 1409 W. Green Street,
Urbana, IL, 61801, USA}
\email{anjid2@illinois.edu}

\address{
Vi Anh Nguyen: Department of Mathematics,
University of Illinois Urbana-Champaign,
Altgeld Hall, 1409 W. Green Street,
Urbana, IL, 61801, USA}
\email{vianhan2@illinois.edu}

\address{
Alexandru Zaharescu: Department of Mathematics,
University of Illinois Urbana-Champaign,
Altgeld Hall, 1409 W. Green Street,
Urbana, IL, 61801, USA and Simion Stoilow Institute of Mathematics of the Romanian Academy, 
P. O. Box 1-764, RO-014700 Bucharest, Romania}
\email{zaharesc@illinois.edu}  
\begin{document}
\setcounter{tocdepth}{1}
\begin{abstract}
We consider the Erd\H{o}s-Moser equation $1^k+2^k+\cdots+(m-1)^k=m^k$ in arithmetic progressions. We prove among other things that when $k=2$, for any solution to exist, the above sum in arithmetic progression must consist of two or four terms. In either case, there are infinitely many solutions that can be completely characterized.    
\end{abstract}
\maketitle

\section{Introduction}\label{sec: Introduction}
The Erd\H{o}s-Moser equation  
\begin{align}
  S_{k}(m) := 1^k + 2^k + \dots + (m-1)^k = m^k,\label{eq:original erdos-moser equation} 
\end{align}
 named after Paul Erd\H{o}s and Leo Moser, has been studied extensively. In \cite{Erdos1949}, Erd\H{o}s conjectured that aside from the trivial solution $1 + 2 = 3$, the equation in \eqref{eq:original erdos-moser equation} has no other solutions. This later became known as the Erd\H{o}s-Moser conjecture. In 1953, Moser \cite{Moser1953} proved that the conjecture is true for odd exponents $k$; for even $k$, he showed that any possible solution $(m, k)$ must satisfy $m > 10^{10^6}$. Moree \cite{Moree2011} later provided an easier proof for Moser's result. 

Moser's bound has since been improved. Butske, Jaje and Mayernik \cite{BJM2000}, by computing certain quantities in Moser’s original proof, showed that $m > 1.485 \times 10^{9 321 155}$ and expressed the hope of reaching the more natural benchmark $10^{10^7}$. Later, Grau and  Oller-Marc\'en in \cite{GrauOllerMarcen2022} defined so-called ``$\mu$-Sondow numbers", a generalization of weak primary pseudoperfect numbers and Giuga
numbers, and showed that Moser's method yields the same bound. Gallot, Moree, and Zudilin \cite{GMZ} substantially improved these bounds, proving that the equation has no solutions for
$m \leq 2.7139 \times 10^{10^{1 667 658 416}}$. This remains the best unconditional result to date.

This problem has been extended and generalized in various ways. Kellner \cite{Kellner2011} studied the generalized Diophantine equation $aS_k(m)=m^k$ and conjectured that there are no solutions with $m > 3$ and integer $a \geq 1$. Meanwhile, Moree \cite{Moree2013} proved that for infinitely many integers $a$, the equation $aS_k(m)=m^k$ has no solution. On the other hand, Sondow and MacMillan \cite{SondowMacMillan2011} reduced the original Erd\H{o}s–Moser equation modulo $k$ and respectively modulo $k^2$ and gave necessary and
sufficient conditions on the solutions to the resulting congruences.
As a corollary to their
results, they obtained a new proof of Moser's result from \cite{Moser1953} for odd exponents $k$. In a later paper \cite{SondowMacMillan2017}, the same authors showed that, under certain conjectures about primary pseudoperfect numbers, no nontrivial solutions exist for $m<10^{10^{20}}$.


In the present paper, we extend the study of Erd\H{o}s-Moser equation to arithmetic progressions. Specifically, we consider the following generalization. 

Fix $q,k\in\N$ and  $b\in\Z$. For any integer $m\geq 2$, consider the equation
\begin{align}\label{eq:generalized E-M equation}
    \sum_{\substack{1 \leq i \leq m -1\\i\equiv b \bmod q}} i^k =m^k.
\end{align}
Unlike the original case of the Erd\H{o}s-Moser equation, by allowing arithmetic progressions, one can have solutions, and in fact infinitely many solutions. 

For instance, when $k = 2$, two solutions are given by $q=7$, $b=5$, and $m=13$, which gives 
\[
5^2+(5+7)^2 = 13^2,
\]
and respectively $q=64$, $b=22$, and $m=276$, which gives
\[
22^2+(22+64)^2+(22+2\times64)^2+(22+3\times64)^2 = 276^2.
\]

Also, for $k = 3$, two solutions are given by $q=808$, $b=317$, and $m=2055$, which gives
\[
317^3+(317+808)^3+(217+2\times808)^3 = 2055^3,
\]
and respectively $q=22215431505$, $b=456326994059$, and $m=2048734872618$, which gives
\begin{align*}
456326994059^3 +(456326994059+22215431505)^3&+(456326994059+2\times22215431505)^3 \\&= 2048734872618^3.
\end{align*}

In what follows, we will mainly be concerned with two types of problems. First, we fix an arithmetic progression. As $m$ goes to infinity, we explore the asymptotic behavior of $k$ in terms of $m$. The proof of our first result uses methods similar to those in \cite[Theorem 1]{GMZ}, which initially inspired our work. Following this direction, we have the following result: 

\begin{thm}\label{thm: main theorem 1}
    Fix integers $q\in\N$ and $b\in\Z$. Suppose $m$ is a positive integer and $k>0$ is a real number satisfying equation \eqref{eq:generalized E-M equation}, and
    assume $b \equiv m - r \bmod q$ for some $r \in \{1,2,\cdots,q\}$. Then we have the asymptotic expansion
    \[
    k = c_1m+c_2+O_q(1/m),
    \]
    where 
    $c_1$ is the solution to $e^{-c_1q}+e^{-c_1r}=1$, and 
    \[
c_2 = - \frac{c_1}{2} \cdot 
\frac{\dfrac{r^2}{1 - e^{-c_1 q}} + \dfrac{2 r q e^{-c_1 q}}{(1 - e^{-c_1 q})^2} + \dfrac{q^2 e^{-c_1 q} (1 + e^{-c_1 q})}{(1 - e^{-c_1 q})^3}}
{\dfrac{r}{1 - e^{-c_1 q}} + \dfrac{q e^{-c_1 q}}{(1 - e^{-c_1 q})^2}}.
\]
\end{thm}
\begin{rem}
 Theorem \ref{thm: main theorem 1} can be proved using the same strategy as in \cite[Theorem 1]{GMZ}. One can define sets $S_n$ and $S_n'$ analogous to those in equations (17) and (18) of \cite{GMZ}, but with the additional restriction that the sum runs over terms in the arithmetic progression $j\equiv m-b\bmod q$. With this newly defined $S_n'$, one can see that \begin{align}
S_n' = \frac{1}{q}\sum_{v=0}^{q-1}\omega^{-rv}A_n(\omega^v z)\Bigg\rvert_{z=e^{-k/m}},\label{eq:Sn' another form}
\end{align}
where $\omega = e^{2\pi i/q}$ and $A_n(t) = (t\frac{d}{dt})^n\frac{t}{1-t}$. Thus, for each $n$,
\[
S_n' -S_n =  \frac{1}{q}\sum_{v=0}^{q-1}\omega^{-rv} B_n(\omega^v z),
\]
where $B_n(t) = (t\frac{d}{dt})^n\frac{t^m}{1-t}$. The rest of the proof proceeds exactly as in \cite[Theorem 1]{GMZ}. However, the computations become more involved for arithmetic progressions. We therefore present a simpler proof of Theorem \ref{thm: main theorem 1}.
\end{rem}
\begin{rem}
    The authors of \cite{GMZ}  obtained an asymptotic formula for $k$ accurate to order $1/m$. We could do the same in Theorem \ref{thm: main theorem 1}, but omit the lower-order terms for simplicity.
\end{rem}

On the other hand, we can fix $k$ in \eqref{eq:generalized E-M equation}, and discuss the possible solutions when allowing $m,b,q$ to vary. Along this direction, taking $k=2$ for simplicity, we obtain the following result.
\begin{thm}\label{thm: main theorem 2}
    When $k=2$, solutions to \eqref{eq:generalized E-M equation} with $m,b,q$ being variables, exist only when the left side consists of two or four terms. Furthermore, in both cases, there are infinitely many solutions to \eqref{eq:generalized E-M equation} and all are ``good" solutions, in the sense that they satisfy
    \[
    \left\lfloor \frac{m-b}{q}\right\rfloor = \left\lfloor\frac{m-b-1}{q}\right\rfloor.
    \]
    Geometrically, each solution corresponds to a point in an infinite set of even-spaced integer points along a line in $\R^3$. See the demonstration in Figure \ref{fig:S_{2,4} stream} below.
\end{thm}
For $k\geq 3$, we make the following remarks. 
\begin{rem} \label{prop:large fixed k}
Let $n$ denote the number of terms on the left side of \eqref{eq:generalized E-M equation}. Assume $k$ is fixed and $m,b,q$ vary. Then, we have the uniform bound
\[
n\leq (k+1)2^{2k+1}.
\]
For sufficiently large $k$, we have 
\[
n<\frac{k}{\log (k+1)}.
\]
Moreover, for $k=3$, there are infinitely many integer solutions to \eqref{eq:generalized E-M equation}. In particular, all such solutions correspond to rational points on certain elliptic curves.  
\end{rem}
\begin{rem}
    Recall that Moser has established that the original Erd\H{o}s-Moser equation in \eqref{eq:original erdos-moser equation} has no solution when the exponent $k$ is odd. By contrast, as Remark \ref{prop:large fixed k} implies, the generalized Erd\H{o}s-Moser equation in arithmetic progressions, defined in \eqref{eq:generalized E-M equation}, has infinitely many solutions when $k=3$.
\end{rem}

\section{Asymptotic behavior of $k$ in terms of $m$ in a fixed arithmetic progression}
\subsection{An Overview of the Erd\H{o}s-Moser equation in Arithmetic Progression with Concrete Examples}
Dividing both sides of \eqref{eq:generalized E-M equation} by $m^k$, we obtain the equivalent form 
\begin{align*}
    \sum_{\substack{i=1\\i\equiv b\bmod q}}^{m-1} \left(\frac{i}{m}\right)^k = 1,
\end{align*}
which can be rewritten as 
\begin{align}
    \sum_{\substack{j=1\\j\equiv m-b\bmod q}}^{m-1} \left(1-\frac{j}{m}\right)^k = 1.\label{eq:generalized E-M equation form 3}
\end{align}
Suppose $b \equiv m - r \bmod q$, or equivalently $m \equiv b + r \bmod q$ for some $r \in \{1,2,\cdots,q\}$. Assume $m = lq + (b + r)$ for some nonnegative integer $l$. Then, the above equation becomes
\begin{align*}
    \left(\frac{b}{m}\right)^k + \left(\frac{b+q}{m}\right)^k + \dots + \left(\frac{b + (l-1)q}{m}\right)^k + \left(\frac{b + lq}{m}\right)^k = \left(\frac{b + lq + r}{m}\right)^k,
\end{align*}
or
\begin{align*}
    \left(\frac{m- (lq+r)}{m}\right)^k + \dots + \left(\frac{m - (2q + r)}{m}\right)^k + \left(\frac{m - (q+r)}{m}\right)^k +  \left(\frac{m - r}{m}\right)^k= 1.
\end{align*}
Rewriting the exponent $k = \frac{k}{m}m$, we obtain 
\begin{align*}
   \left(\left(\frac{m- (lq+r)}{m}\right)^m\right)^{k/m} + \dots &+ \left(\left(\frac{m - (2q + r)}{m}\right)^m\right)^{k/m}\\ 
   &+ \left(\left(\frac{m - (q + r)}{m}\right)^m\right)^{k/m} + \left(\left(\frac{m -  r}{m}\right)^m\right)^{k/m} = 1.
\end{align*}
Using L'H\^{o}pital's rule, one can show that 
\[
\displaystyle{\lim_{m \to \infty} \left(\frac{m-r}{m}\right)^m} = \frac{1}{e^r}.
\]
  We have $\displaystyle{\lim_{m \to \infty} m \ln \left(\frac{m-r}{m}\right)} = \displaystyle{\lim_{m \to \infty} \frac{\ln \left(\frac{m-r}{m}\right)}{m^{-1}}} = \displaystyle{\lim_{m \to \infty} \frac{\frac{m}{m-r}rm^{-2}}{-m^{-2}}} = \displaystyle{\lim_{m \to \infty} \frac{-rm}{m - r}} = -r$. Therefore, the result follows. Similarly, one can show that
  \[
  \displaystyle{\lim_{m \to \infty} \left(\frac{m-(gq+r)}{m}\right)^m} = \frac{1}{e^{gq + r}}
  \]
   for any $g=1,2,\cdots,l$. Thus, as $k, m \to \infty$, we achieve
\begin{align*}
\left(\frac{1}{e^r}\right)^{k/m}
+ \left(\frac{1}{e^{q+r}}\right)^{k/m}
+ \left(\frac{1}{e^{2q+r}}\right)^{k/m} + \dots = 1.
\end{align*}
Let $u = \displaystyle{\lim_{k, m \to \infty} e^{-k/m}}$. Note that $u$ is always nonnegative. Thus,
\begin{align*}
    u^r + u^{q+r} + u^{2q+r} + \dots = 1, 
\end{align*}
or
\begin{align*}
    u^r(1 + u^q + u^{2q} + \dots) = \frac{u^r}{1 - u^q} = 1.
\end{align*}
The problem therefore reduces to finding solutions of the polynomial
\begin{align}
    u^q + u^r - 1 = 0. \label{eq:solution of a polynomial}
\end{align}
Let $f(x):=x^q + x^r - 1$. Then $f(0)=-1$ and $f(1)=1$, so there exists a solution $u\in (0,1)$ satisfying $f(u)=0$. Moreover, since $f'(x) = qx^{q-1}+rx^{r-1}>0$ for $x>0$, this solution $u$ is unique. Note that \eqref{eq:solution of a polynomial} is solvable by radicals when $q<5$. When $q\geq 5$, this might not be solvable by radicals. For instance, when $q=5$ and $r=1$, the Galois group of this polynomial in \eqref{eq:solution of a polynomial} is $S_5$, which is not solvable.

We now discuss some special cases. When $r=1$, \eqref{eq:solution of a polynomial} becomes 
\begin{align}
    u^q + u - 1 = 0. \label{eq:solution of a polynomial r=1}
\end{align}
When $q=1$, we have $u=1/2$, which reduces to the case in \cite[Theorem 1]{GMZ}. When $q=2$, we will have $u=\frac{\sqrt{5}-1}{2}$, which is the inverse Golden ratio. When $q\to\infty$, we have $u^q\to 0$, which implies $u\to 1$. Assume $u=1-\delta$, then \eqref{eq:solution of a polynomial r=1} becomes
\[
(1-\delta)^q=\delta.
\]
Taking logarithms on both sides, we have
\[
q\log(1-\delta) = \log \delta.
\]
Let $\delta = \dfrac{1}{M}$. The equation becomes
\[
-q\log \left(1- \dfrac{1}{M}\right) = \log M.
\]
Using the Taylor expansion of $-\log(1-x)$, we have
\[
q\left(\dfrac{1}{M} + \dfrac{1}{2M^2} + \dfrac{1}{3M^3} + \dots \right) = \log M,
\]
or
\[
q\left(1 + \dfrac{1}{2M} + \dfrac{1}{3M^2} + \dots \right) = M\log M.
\]
This yields
\begin{align}
q\left(1+O\left(\dfrac{1}{M}\right)\right) = M\log M. \label{M in terms of q}
\end{align}
Again, taking logarithms of both sides, we have
\begin{align*}
\log q + O\left(\dfrac{1}{M}\right) = \log M + \log \log M.
\end{align*}
Thus, we can write
\begin{align}
\log M = \log q + O(\log \log M).\label{log M in terms of log q}
\end{align}
Using \eqref{M in terms of q} and \eqref{log M in terms of log q}, we have
\begin{align*}
    M = \dfrac{q\left(1 + O\left(\dfrac{1}{M}\right)\right)}{\log q\left(1 + O\left(\dfrac{\log \log M}{\log q}\right)\right)}.
\end{align*}
Finally, using $\delta = \dfrac{1}{M}$ and noting that $M>\log q$, we have asymptotically
\begin{align*}
    \delta = \dfrac{\log q}{q}\left(1 + O\left(\dfrac{\log \log \log q }{\log q}\right)\right).
\end{align*}
Thus, when $q\to\infty$, we have $u\to 1-\log q/q$.

When $r=q$, \eqref{eq:solution of a polynomial} becomes
\[
2u^q = 1,
\]
which implies that $u=2^{-1/q}$. In this case, when $q=1$, we have $u=1/2$, which reduces to the case in \cite{GMZ}; when $q=2$, $u=\sqrt{2}/2$, and finally, when $q\to\infty$, we have $u\to 1^{-}$.

\subsection{Proof of Theorem \ref{thm: main theorem 1}} Now we will prove Theorem \ref{thm: main theorem 1}. Note that
\begin{align}
    \left(1-\frac{j}{m}\right)^k &= e^{k\log(1-\frac{j}{m})}\notag\\
    &= e^{k(-\frac{j}{m}-\frac{j^2}{2m^2}-\frac{j^3}{3m^3}+O(j^4m^{-4}))}, \label{eq: expansion of (1-j/m) kth power}
\end{align}
where the second and third equalities follow from the Taylor expansion of $\log(1-j/m)$.

Krzysztofek \cite{Krzysztofek1966} shows that $k\gg m$ holds even without arithmetic progression, so the same bound applies in our arithmetic progression setting. Applying \cite[ Lemma 1]{GMZ}, when $m$ goes to infinity, we have
\begin{align}
    e^{-\frac{j^2k}{2m^2}+O(kj^3m^{-3}))} = 1-\frac{j^2k}{2m^2}+O(kj^3m^{-3}).\label{eq:equality of exp and lemma 1}
\end{align}
Using Taylor expansion on the left-hand side, we see that \eqref{eq:equality of exp and lemma 1} is valid only if the leading term of the squared argument on the left side is absorbed by the error term on the right side, i.e.
\[
\frac{j^4k^2}{m^4} = O\left(\frac{kj^3}{m^3}\right),
\]
which simplifies to $jk = O(m)$, forcing $k=O(m)$. Thus, $k\sim c_1m$ for some constant $c_1\neq 0$.

We make the ansatz $k=c_1 m+c_2+O(1/m)$. Then, 
\[
e^{k(-\frac{j}{m}-\frac{j^2}{2m^2}+O(m^{-3}))} = e^{-c_1 j + \frac{-c_1 j^2 - 2 c_2 j}{2 m}  + O\left(m^{-2}\right)},
\]
and thus,
\[
 \sum_{\substack{j=r\\j\equiv r\bmod q}}^{lq+r}e^{-c_1 j + \frac{-c_1 j^2 - 2 c_2 j}{2 m}  + O\left(m^{-2}\right)} =  1.
\]
This is equivalent to
\begin{align*}
    &\sum_{\substack{j=r\\j\equiv r\bmod q}}^{lq+r}e^{-c_1j} -\left( \sum_{\substack{j=r\\j\equiv r\bmod q}}^{lq+r}e^{-c_1 j}\frac{c_1 j^2}{2m}+\sum_{\substack{j=r\\j\equiv r\bmod q}}^{lq+r}e^{-c_1 j}\frac{c_2 j}{m}\right)+O(m^{-2}) = 1.
\end{align*}
If we write $j=r+wq$, then this is equivalent to
\begin{align}
    e^{-c_1r}\sum_{w=0}^{l} e^{ -c_1w q}   -e^{-c_1r}\sum_{w=0}^{l} e^{ -c_1w q}\left( \frac{c_1 (r + w q)^2}{2m} + \frac{c_2 (r + w q)}{m}+O\Big(\frac{1}{m^2}\Big)\right)    = 1.\label{eq:equality before finding c1 c2}
\end{align}
When $l$ tends to infinity, 
\[
\lim_{l\to\infty}\sum_{w=0}^{l} e^{ -c_1w q} = \frac{1}{1-e^{-c_1q}}.
\]
Since the main term equals to $1$, we achieve that 
\[
\frac{e^{-c_1r}}{1-e^{-c_1q}}=1,
\]
so $c_1$ is the solution to $e^{-c_1q}+e^{-c_1r}=1$, and 
\eqref{eq:equality before finding c1 c2} gives us
\begin{align}
    -e^{-c_1r}\sum_{w=0}^{l} e^{ -c_1w q}\left( \frac{c_1 (r + w q)^2}{2m} + \frac{c_2 (r + w q)}{m}\right)  = O\Big(\frac{1}{m^2}\Big).\label{eq:equality before finding c2}
\end{align}
Recall that $\sum_{w=0}^\infty wx^w = \frac{x}{(1-x)^2}$ and $\sum_{w=0}^{\infty} w^2 x^w = \frac{x(1+x)}{(1-x)^3}$. Therefore, for the second-order terms, we have
\begin{align*}
    \lim_{l\to\infty}\sum_{w=0}^{l} e^{ -c_1w q} (r+wq)=\frac{r}{1 - e^{-c_1 q}} + \frac{q e^{-c_1 q}}{(1 - e^{-c_1 q})^2}
\end{align*}
and
\begin{align*}
   \lim_{l\to\infty} \sum_{w=0}^{l} e^{ -c_1w q} \frac{(r+wq)^2}{2} 
= \frac{r^2}{1 - e^{-c_1 q}} + \frac{2 r q e^{-c_1 q}}{(1 - e^{-c_1 q})^2} + \frac{q^2 e^{-c_1 q} (1 + e^{-c_1 q})}{(1 - e^{-c_1 q})^3}.
\end{align*}
Since the left-hand side of \eqref{eq:equality before finding c2} must be $O(m^{-2})$, we require
\[
\sum_{w=0}^{\infty} \left( \frac{c_1}{2} (r + w q)^2 + c_2 (r + w q) \right) e^{-c_1 q w} = 0,
\]
and thus, 
\[
c_2 = - \frac{c_1}{2} \cdot 
\frac{\dfrac{r^2}{1 - e^{-c_1 q}} + \dfrac{2 r q e^{-c_1 q}}{(1 - e^{-c_1 q})^2} + \dfrac{q^2 e^{-c_1 q} (1 + e^{-c_1 q})}{(1 - e^{-c_1 q})^3}}
{\dfrac{r}{1 - e^{-c_1 q}} + \dfrac{q e^{-c_1 q}}{(1 - e^{-c_1 q})^2}}.
\]
This confirms our ansatz, so indeed $k = c_1m+c_2+O(1/m)$. 
In the particular case, when $q=r=1$, $b=0$ and $l = m-1$, we have 
\[
c_1 = \log 2, \quad c_2 = -\frac{3}{2}\log 2,
\]
which matches the result in \cite[Theorem 1]{GMZ}.

\section{Erd\H{o}s-Moser Equation with a fixed exponent $k$}
We now consider the direction when we have a fixed exponent $k$. Let $n\geq 2$ denote the number of terms on the left-hand side of \eqref{eq:generalized E-M equation}. Note that \eqref{eq:generalized E-M equation} is equivalent to 
\begin{align}
    b^k+(b+q)^k+\cdots+\left(b+\left\lfloor\frac{m-1-b}{q}\right\rfloor q\right)^k = m^k,\label{eq:fixed number of terms with k}
\end{align}
with $1\leq b\leq q$ and $\lfloor (m-1-b)/q\rfloor=n$. For the rest of the section, we will fix $k=2$, and view $m,b,q$ as variables by letting $x=q,y=b$, and $z=m$. 

The problem thus reduces to finding integer solutions of
\begin{align}\label{0.3}
    y^2+(y+x)^2+\cdots+(y+\lfloor\frac{z-1-y}{x}\rfloor x)^2 - z^2 = 0,
\end{align}
where $1\leq y\leq x$ and $z>0$. 

Let $S_2$ be the surface defined by \eqref{0.3}. We decompose $S_2$ as a disjoint union of surfaces $S_{2,n}$, i.e., 
\begin{align}
    S_2  = \bigcup_{n\geq 2} S_{2,n},\label{eq:definition of S2}
\end{align}
where 
\begin{align}
    S_{2,n} &= \{(x,y,z)\in\R^3:  y^2+(y+x)^2+\cdots+(y+(n-1) x)^2 - z^2 = 0; \notag\\
    &\quad1\leq y\leq x,(n-1)x+y+1\leq z< nx+y+1\}.\label{eq:S(2,n) definition}
\end{align}

Let $\mathcal{D}_{2,n}$ be the projection of $S_{2,n}$ onto the $(x,y)$ plane, i.e.
\begin{align}\label{D(k,n)}
   \mathcal{D}_{2,n} &:= \{(x,y)\in\R^2: 1\leq y\leq x,\notag\\
   &\quad((n-1)x+y+1)^2\leq y^2+(y+x)^2+\cdots+(y+(n-1) x)^2 < (nx+y+1)^2\}.
\end{align}

\subsection{Parametric Representation}
Using the formula for square pyramidal numbers, we obtain the following equality:
\begin{align*}
    y^2+(y+x)^2+\cdots+(y+(n-1) x)^2  = ny^2+\frac{n(n-1)(2n-1)}{6}x^2+n(n-1)xy.
\end{align*}
Thus, $\mathcal{D}_{2,n}$ and $S_{2,n}$ simplify to
\begin{align*}
   \mathcal{D}_{2,n} = &\{(x,y)\in\R^2: 1\leq y\leq x,\\
   &\quad((n-1)x+y+1)^2\leq ny^2+\frac{n(n-1)(2n-1)}{6}x^2+n(n-1)xy< (nx+y+1)^2\},
\end{align*}
and \begin{align*}
   S_{2,n} &= \{(x,y,z)\in\R^3: z^2=ny^2+\frac{n(n-1)(2n-1)}{6}x^2+n(n-1)xy, 1\leq y\leq x,\\
   &\quad(n-1)x+y+1\leq z
   < nx+y+1\}.
\end{align*}

Note that we can write the real solutions in $S_{2,n}$ in parametric representations. For instance, when $n=2$, i.e. there are only two terms in the arithmetic progression, 
the parametric representation of $S_{2,2}$ is
\begin{align*}
\begin{cases}
   x &= u - v, \\
    y &= v, \\
    z &= \sqrt{u^2 + v^2},
\end{cases} 
\end{align*}
where $ 1 \leq v \leq \frac{u}{2}$ and $u + 1 \leq \sqrt{u^2 + v^2} < 2u-v + 1$. This can be rewritten as
\begin{align*}
\begin{cases}
   x &= u - 1-t(\frac{u}{2}-1), \\
    y &= 1+t(\frac{u}{2}-1), \\
    z &= \sqrt{u^2 + (1+t(\frac{u}{2}-1))^2},
\end{cases} 
\end{align*}
where $0\leq t\leq 1$, $u\geq 2$, $\sqrt{2u+1}\leq 1+t(\frac{u}{2}-1)$ and $3u^2 - 2u^2 t + 3ut + 2t + 3>0$. For example, taking $t=0.8$ and $u=12$ yields $x=7,y=5$, and $z=13$.

\subsection{Exploration on the Solution Sets in $S_{2,n}$}
We first explore $\mathcal{D}_{2,n}$ and $S_{2,n}$ graphically to gain insight. 
\begin{figure}[ht]
\centering
\includegraphics[width=10cm, height=8cm]{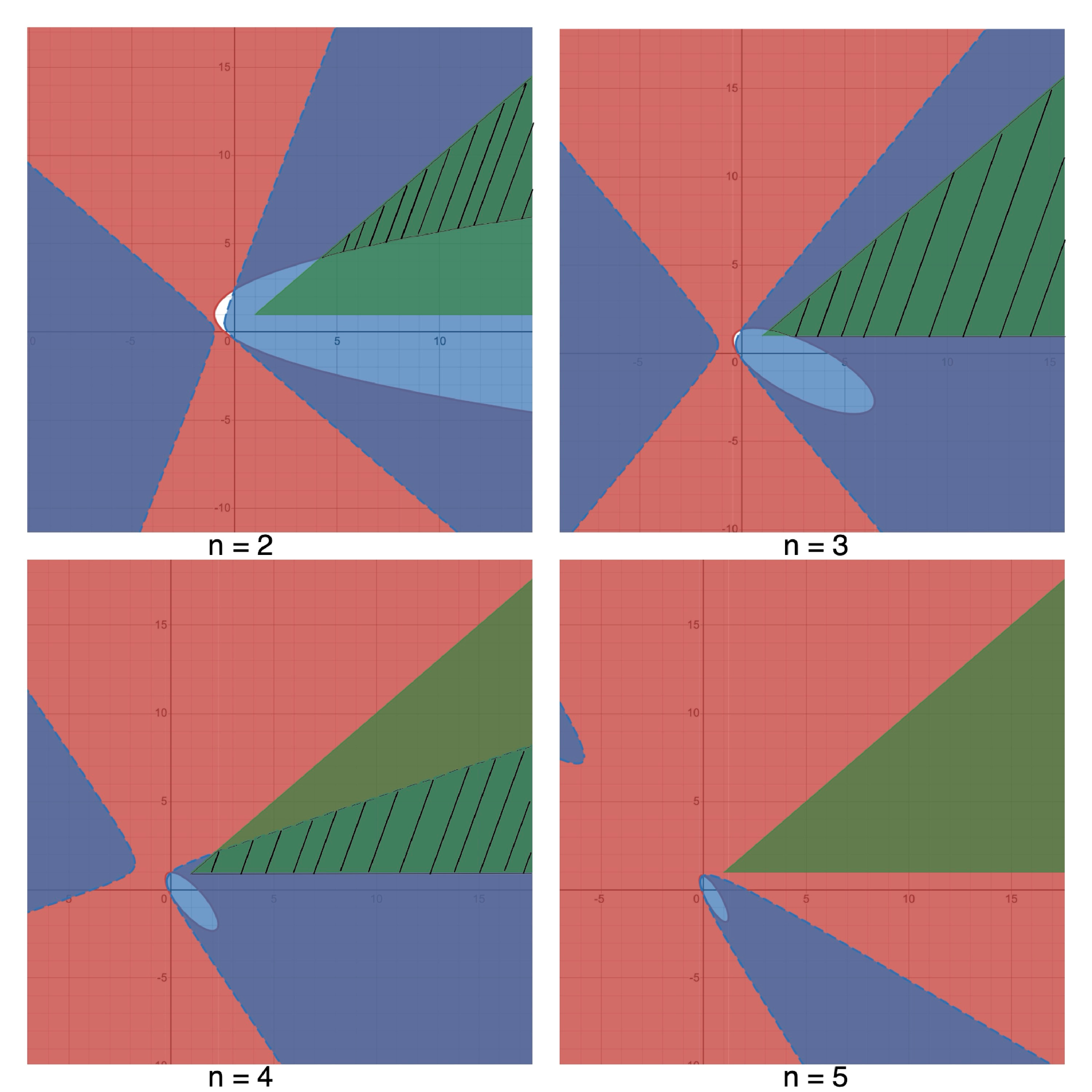} 
\caption{Examples of $\mathcal{D}_{2,n}$ (the part that is shaded with black lines)}\label{fig:D(2,n)}
\end{figure}
\begin{figure}[ht]
\centering
\includegraphics[width=15cm, height=6cm]{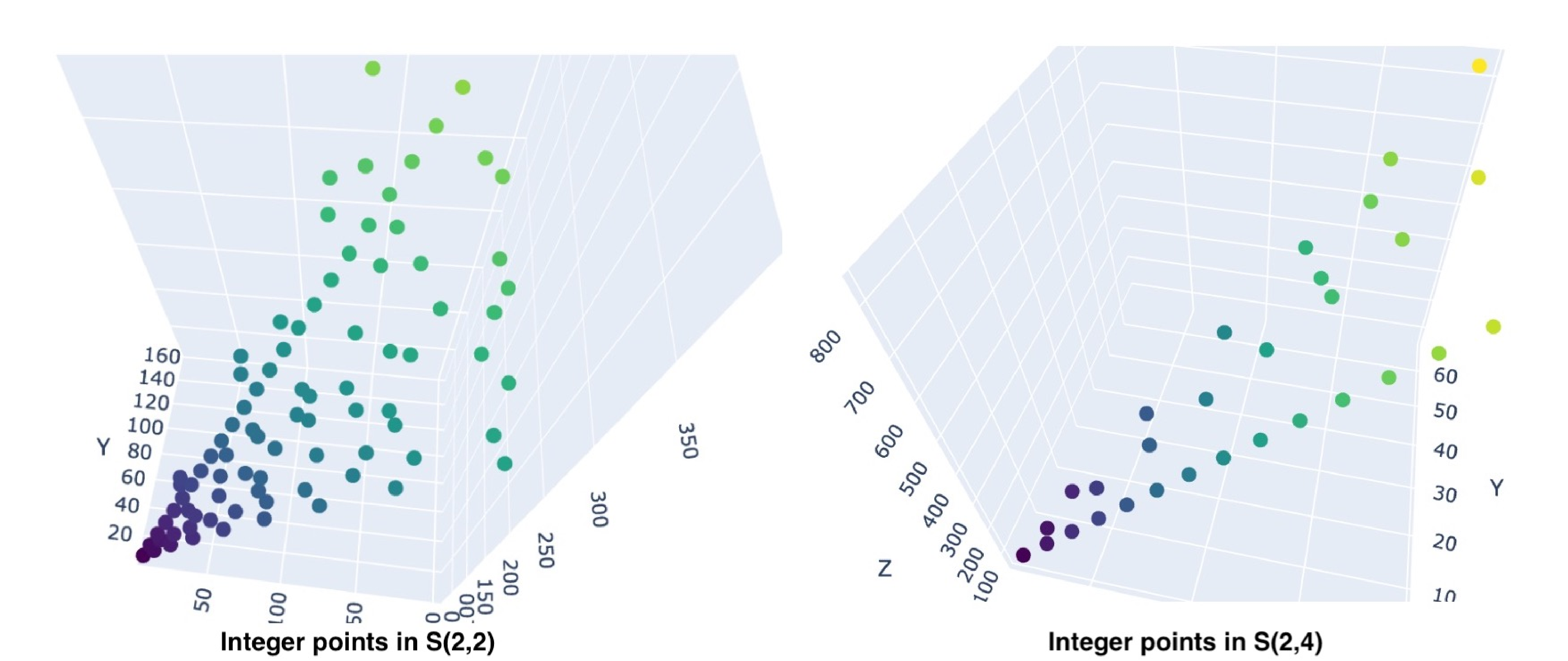} 
\caption{Examples of $S_{2,n}$}\label{fig:S_{2,n}}
\end{figure}
 Figure \ref{fig:D(2,n)} shows that for $n=2,3,4$, there are infinitely many integer points $(x,y)\in\N^2$ in $\mathcal{D}_{2,n}$, but for $n\geq 5$, no solutions exist. Correspondingly, one might expect infinitely many integer solutions in $S_{2,n}$ for $n=2,3,4$. This is confirmed true for $n=2$ and $n=4$, but checking the first $20,000\times 20,000$ lattice points in $\mathcal{D}_{2,3}$ reveals no integer points in $S_{2,3}$. The integer points in $S_{2,n}$ for $n=2,4$ are shown in Figure \ref{fig:S_{2,n}}.

Based on the above observations, we propose the following conjectures.
\begin{conj}\label{conjecture 1}
  $\mathcal{D}_{2,n}$ is nonempty only when $n=2,3,4$. 
\end{conj}

\begin{conj}\label{conjecture 2}
    There does not exist any integer point in $S_{2,3}$. 
\end{conj}

It's easy to show that $S_{2,2}$ and $S_{2,4}$ contain infinitely many integer points. Suppose $(x,y,z)$ is an integer solution in $S_{2,n}$ for $n=2$ or $4$. Observe that if for $j\in\N$, the integer point $(x,y,z)$ satisfies 
\[
\left\lfloor \frac{z-y-1/j}{x}\rfloor = \lfloor\frac{z-y-1}{x}\right\rfloor,
\]
then $(jx,jy,jz)$ is also a solution to $S_{2,n}$. Consequently, if 
\begin{align}
\left\lfloor \frac{z-y}{x}\right\rfloor = \left\lfloor\frac{z-y-1}{x}\right\rfloor,\label{eq:condition for good solution}
\end{align}
then one can generate infinitely many solutions of the form $(jx,jy,jz)$ for any $j\in\N$. These solutions correspond to evenly spaced lattice points along a line, as shown in Figure \ref{fig:S_{2,4} stream} below. We call such a solution a ``good" solution. Moreover, we call any ``good" solution $(x,y,z)$ with gcd$(x,y,z)=1$ a ``minimally good" solution. 

For $S_{2,2}$, note that $x=7,y=5,z=13$ is a ``minimally good" solution. Indeed, for any positive integer $j$, we have $\lfloor \frac{z-y-1/j}{x}\rfloor = \lfloor \frac{8-1/j}{7}\rfloor = \lfloor\frac{8-1}{7}\rfloor$, so $(x,y,z)=(7j,5j,13j)$ is an integer solution in $S_{2,2}$. Similarly, $x=14,y=1,z=54$ is a ``minimally good" solution in $(2,4)$. Therefore, both $S_{2,2}$ and $S_{2,4}$ contain infinitely many integer points. 

In Figure \ref{fig:S_{2,4} stream}, all marked points are ``minimally good" solutions in $S_{2,4}$, including the two circled points. Checking the first $20000\times20000$ lattice points in the $xy$-plane implies that all integer solutions in $S_{2,2}$ and $S_{2,4}$ are ``good" solutions. This leads to the following conjecture.
\begin{conj}\label{conjecture 3}
    All integer points in $S_{2,2}$ and $S_{2,4}$ are ``good" solutions. 
\end{conj}

\begin{figure}[ht]
\centering
\includegraphics[width=10cm, height=8cm]{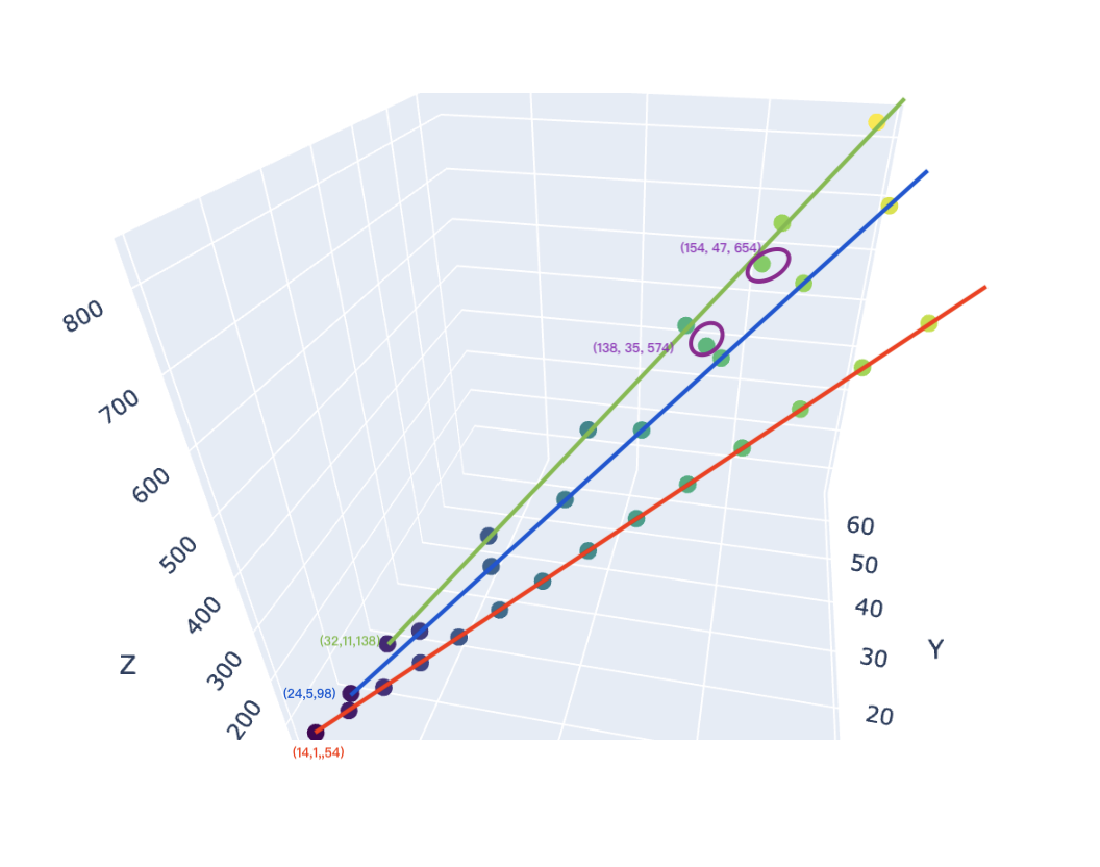} 
\caption{Integer solutions in $S_{2,4}$ }\label{fig:S_{2,4} stream}
\end{figure}

\section{Proof of Theorem \ref{thm: main theorem 2}}
To prove Theorem \ref{thm: main theorem 2}, it suffices to prove Conjectures \ref{conjecture 1}--\ref{conjecture 3}.
\subsection{Proof of Conjecture \ref{conjecture 1}}
It's easy to see that there are infinitely many integer points in $\mathcal{D}_{2,n}$ for $n=2,3,4$. In fact, the diagonal points $(t,t)\in\Z^2$ with $t\geq 5$ lie in both $\mathcal{D}_{2,2}$ and $\mathcal{D}_{2,3}$, while the points $(t,1)$ for $t\in\N$ lie in $\mathcal{D}_{2,4}$. For $n\geq 5$, define
\[
L:=((n-1)x+y+1)^2,\quad
M:=ny^2+A x^2+n(n-1)xy,\quad\textrm{and}\quad
U:=(nx+y+1)^2,
\]
where
\[
A=\frac{n(n-1)(2n-1)}{6}.
\]
Then $(x,y)\in\mathcal{D}_{2,n}$ if and only if $1\leq y\leq x$ and $L\leq M<U$. To prove that $\mathcal{D}_{2,n}$ has no integer points, it suffices to show that $U\leq M$. To do so, define 
\[f(x,y) := U-M = (1 - n)y^2  + ((3n - n^2)x+2) y +(n^2 - A)x^2+ 2nx  + 1.
\]
For fixed $x$ and $n\geq 2$, the quadratic function $f(x,y)$ is concave  in $y$ on $[1,x]$. The vertex $v$ of $f(y)$ occurs at 
\[
v = -\frac{n(3-n)x+2}{2(1-n)}=\frac{n(3-n)x+2}{2(n-1)},
\]
which is negative for $n\geq 5$. Therefore, $f(y)$ achieves its maximum value at $y=1$ for any fixed $x$. However,
\[
F(x,1)=\frac{n(13-2n^2+3n)}{6}x^2+n(3-n)x-(n-1)+2nx+3,
\]
which is negative for all integers $x\geq 1$ when $n\geq 5$. Thus $U\leq M$, as desired. This completes the proof of Conjecture \ref{conjecture 1}. 

\subsection{Proof of Conjecture \ref{conjecture 2}}
Integer points in $S_{2,3}$ exactly satisfy 
\begin{enumerate}
    \item $5x^2+6xy+3y^2 = z^2$
    \item $1\leq y\leq x$
    \item $2x+y+1\leq z< 3x+y+1$.
\end{enumerate}
Reducing the first equation modulo $3$, we obtain
\[
z^2 \equiv 2x^2 \bmod 3.
\]
Note that $2x^2$ is a quadratic residue modulo $3$ if and only if $x\equiv 0 \bmod 3$. Therefore, write $x=3t$. Substituting into the first equation gives that
\begin{align}
z^2 = 5(3t)^2 + 6(3t)y + 3y^2 = 3(15t^2+6ty+y^2).\label{eq:after substitute x=3t}
\end{align}
Thus $3\mid z$, so the right-hand side of \eqref{eq:after substitute x=3t} is divisible by $9$, forcing $3\mid y$. Consequently, any integer point $(x,y,z)$ in $S_{2,3}$ must be of the form $(3x_1,3y_1,3z_1)$ with $x_1,y_1,z_1\in\N$. Substituting back into the first equation in the conditions, we find that $(x_1,y_1)$ satisfies
\begin{enumerate}
    \item $5x_1^2+6x_1y_1+3y_1^2=z_1^2$
    \item $\frac{1}{3}\leq y_1\leq x_1$
    \item $2x_1+y_1+\frac{1}{3}\leq z_1< 3x+y+\frac{1}{3}$.
\end{enumerate}
Since $x_1,y_1,z_1\in\N$, it follows that $(x_1,y_1,z_1)$ is also in $S_{2,3}$. Iterating this argument produces an infinite strictly descending sequence of positive integer solutions
\[
(x,y,z) \mapsto (x_1,y_1,z_1) \mapsto (x_2,y_2,z_2) \mapsto \cdots,
\]
which is impossible. Hence, $S_{2,3}$ contains no integer solutions. This completes the proof of Conjecture \ref{conjecture 2}.

\subsection{Proof of Conjecture \ref{conjecture 3}}
Suppose $(x,y,z)$ is an integer point in $S_{2,n}$ for $n=2$ or $n=4$ such that $(x,y,z)$ is not ``good". Then, \eqref{eq:condition for good solution} fails. This implies that there exists an integer $u$ such that 
\[\frac{z-y-1}{x} < u \leq \frac{z-y}{x}\quad\textrm{or}\quad \frac{z-y-1}{x} \leq u <\frac{z-y}{x}.
\]
Without loss of generality, we consider the first scenario. This integer $u$ is unique and equals $\lfloor(z-y)/x\rfloor$. Therefore, we have
\[
\frac{z-y-1}{x} < \left\lfloor \frac{z-y}{x} \right\rfloor. 
\]
Since
$\left\{\frac{z-y}{x}\right\} = \frac{z-y}{x} - \left\lfloor \frac{z-y}{x} \right\rfloor$, we arrive at
\[\left\{\frac{z-y}{x}\right\} < \frac{1}{x},
\]
which implies that $x\mid (z-y)$. Thus, it suffices to prove that for any integer point in $S_{2,2}$ or $S_{2,4}$, we have $x\nmid (z-y)$.

When $n=2$, the corresponding equation is 
\begin{align}
z^2 = x^2+2xy+2y^2,\label{eq: S2,2 corresponding equation}
\end{align}
where $1\leq y\leq x$ and $x+y+1\leq z < 2x+y+1$. Now suppose $\kappa = (z-y)/x \in \Z$. From the inequality $x+y+1\leq z< 2x+y+1$, we have 
\[
x+1 \leq z-y < 2x+1 ,
\]
which implies 
\[
1+\frac{1}{x}\leq \kappa < 2+\frac{1}{x}.
\]
Since $\kappa\in\Z$, we must have $\kappa = 2$, and hence $z = y+2x$. Substituting into \eqref{eq: S2,2 corresponding equation} and simplifying yields
\[
x^2-2xy-2y^2=0.
\]
Viewing this as a polynomial in $x$, the solution is 
\[
x = y(1\pm \sqrt{3}),
\]
which cannot be an integer for any positive integer $y$. This contradiction shows that $x\nmid (z-y)$ for integer points in $S_{2,2}$.

When $n=4$, the equation is
\begin{align}
z^2 = 4y^2+14x^2+12xy,\label{eq: S2,4 corresponding equation}
\end{align}
with the inequalities $1\leq y\leq x$ and $3x+y+1\le z < 4x+y+1$. Following the same strategy, we obtain $z=y+4x$. Substituting into \eqref{eq: S2,4 corresponding equation} gives
\[
2x^2-4xy-3y^2=0.
\]
Solving for $x$ in terms of $y$ yields $x = y\big(\frac{1\pm\sqrt{10}}{2}\big)$, which is again impossible for positive integer $y$, so $x\nmid (z-y)$. This finishes the proof of Conjecture \ref{conjecture 3}. Combining the results in Conjectures \ref{conjecture 1}--\ref{conjecture 3}, we complete the proof of Theorem \ref{thm: main theorem 2}.

\section{Some remarks when $k\geq 3$}
\subsection{An upper bound on $n$ depending on $k$.}
Let $x=q,y=b,z=m$, and $n=\lfloor (m-1-b)/q\rfloor$ as in  \eqref{eq:fixed number of terms with k}. We define the surface $S_{k}$, analogous to $S_{2}$ in \eqref{eq:definition of S2} and \eqref{eq:S(2,n) definition}, by
\begin{align*}
    S_k = \bigcup_{n\geq 2} S_{k,n},
\end{align*}
where
\begin{align}
    S_{k,n} &= \{(x,y,z)\in\R^3:  y^k+(y+x)^k+\cdots+(y+(n-1) x)^k  = z^k; \notag\\
    &\quad1\leq y\leq x,(n-1)x+y+1\leq z< nx+y+1\}.
\end{align}
Then we have
  
\begin{align}
    \sum_{j=0}^{n-1} (y+jx)^k\geq \sum_{j=0}^{n-1}(jx)^k\geq x^k\int_{0}^{n-1}t^kdt = \frac{x^k(n-1)^{k+1}}{k+1}.\label{eq: inequality from below}
\end{align}
On the other hand, the assumption that $z<nx+y+1\leq (n+1)x+1$ implies that 
\begin{align}
    z^k<((n+1)x+1)^k\leq (n+2)^kx^k.\label{eq:inequality from above}
\end{align}
Combining \eqref{eq: inequality from below} and \eqref{eq:inequality from above}, we obtain 
\begin{align*}
    \frac{(n-1)^{k+1}}{k+1}\leq (n+2)^k.
\end{align*}
Using the fact that $(n-1)^{k+1}\geq (n/2)^{k+1}$ for $n\geq 2$ and $(n+2)^k \leq (2n)^k$, we obtain
\begin{align*}
   \frac{(n/2)^{k+1}}{k+1}\leq (2n)^k,
\end{align*}
which implies that 
\begin{align}
n\leq (k+1)2^{2k+1}.\label{eq:uniform upper bound}    
\end{align}
Therefore, for any fixed $k$, any integer solution to $S_k$ must have $n$ bounded. For large $k$, using 
\begin{align}
    \sum_{j=0}^{n-1} (y+jx)^k\leq \sum_{j=0}^{n-1}((j+1)x)^k\leq x^k\int_{1}^{n}t^kdt \leq \frac{x^kn
    ^{k+1}}{k+1}\label{eq: inequality from above 2}
\end{align}
and the inequality $(n-1)x+y+1\leq z$, we obtain
\[
k\log(1-1/n)<-\log(k+1).
\]
For large $n$, the left side is asymptotic to $-k/n$, yielding the bound 
\[
n<\frac{k}{\log(k+1)}. 
\]

\subsection{A discussion on $k=3$.} For $k=3$, we have a tighter upper bound on $n$: 
\[
(n-1)^4\leq 4(n+2)^3,
\]
so it suffices to check for $n=2,3,\cdots,10$. The case $n=2$ is eliminated by Fermat's Last Theorem. For the remaining $n$, numerically checking with $x,y\leq 5000$ reveals that when $n=3$, integer solutions of the form $(808j,317j,2055j)$ exist for $j\in\N$. We therefore explore this case further.

We have \begin{align}
    S_{3,3} &= \{(x,y,z)\in\R^3:  y^3+(y+x)^3+(y+2x)^3  = z^3; \notag\\
    &\quad1\leq y\leq x,2x+y+1\leq z< 3x+y+1\}.\label{eq: S(3,3)}
\end{align}
Expand the curve equation in \eqref{eq: S(3,3)} and let $a=x+y$ and $b=x$. Then the equation is equivalent to
\begin{align}
z^3 = 3a^3+6ab^2.\label{eq:S(3,3) equivalent form}
\end{align}
 Divide both sides by $a^3/216$, and set $u=z/a$ and $v=b/a$. We then obtain
\[
(36v)^2 =  (6u)^3-648.
\]
Now take $W = 36v$ and $R=6u$, we arrive at
\[
W^2 = R^3-648,
\]
which is an elliptic curve, since it has nonzero discriminant. Moreover, the conditions $1\leq y\leq x$ and $2x+y+1\leq z<3x+y+1$ imply that $18\leq W<36$ and $6+W/6\leq R<6+W/3$. Denote this elliptic curve by $E$. Then, any integer point in $S_{3,3}$ must be a rational point on $E$. Using SageMath, one can compute that the Mordell-Weil group $E(\Q)$ is generated by $P = E(9,9)$ and $Q=E(18,72)$. Therefore, any rational point on $E$ must be of the form
\begin{align}
(R,W) = n_1P+n_2Q,\label{eq:generators of the elliptic curve}
\end{align}
with $n_1,n_2\in\Z$.

The solution $(x,y,z)=(808,317,2055)$ we found earlier corresponds to the rational point $(R,W) = (\frac{274}{25},\frac{3232}{125})$ on $E$, and we have
\[
(\frac{274}{25},\frac{3232}{125}) = -2P+Q.
\]
 Using \eqref{eq:generators of the elliptic curve}, we can find other integer solutions to $S_{3,3}$ similarly. For instance, when $n_1=-1$ and $n_2=-4$, we obtain 
$(x,y,z)=(722215431505,456326994059,2048734872618)$.

\printbibliography

@article {GMZ,
    AUTHOR = {Gallot, Y. and Moree, P. and Zudilin, W.},
     TITLE = {The {E}rd\H os-{M}oser equation {$1^k+2^k+\dots+(m-1)^k=m^k$}
              revisited using continued fractions},
   JOURNAL = {Math. Comp.},
  FJOURNAL = {Mathematics of Computation},
    VOLUME = {80},
      YEAR = {2011},
    NUMBER = {274},
     PAGES = {1221--1237},
      ISSN = {0025-5718,1088-6842},
   MRCLASS = {11D61 (11A55 11Y65)},
  MRNUMBER = {2772120},
MRREVIEWER = {A.\ Peth\H o},
       DOI = {10.1090/S0025-5718-2010-02439-1},
       URL = {https://doi-org.proxy2.library.illinois.edu/10.1090/S0025-5718-2010-02439-1},
}

@article {Krzysztofek1966,
    AUTHOR = {Krzysztofek, B.},
     TITLE = {The equation {$1\sp{n}+2\sp{n}+\cdots
              +m\sp{n}=(m+1)\sp{n}\cdot k$}},
   JOURNAL = {Wy\.z. Szko\l. Ped. w Katowicach---Zeszyty Nauk. Sekc. Mat.},
  FJOURNAL = {Wy\.zsza Szko\l a{} Pedagogiczna w Katowicach. Zeszyty
              Naukowe. Sekcja Matematyki},
    NUMBER = {5},
      YEAR = {1966},
     PAGES = {47--54},
      ISSN = {2082-9566},
   MRCLASS = {10B99},
  MRNUMBER = {364097},
MRREVIEWER = {B.\ Nov\'ak},
}

@article{Kellner2011,
    AUTHOR = {Kellner, B. C.},
     TITLE = {On stronger conjectures that imply the {E}rd\H os-{M}oser
              conjecture},
   JOURNAL = {J. Number Theory},
  FJOURNAL = {Journal of Number Theory},
    VOLUME = {131},
      YEAR = {2011},
    NUMBER = {6},
     PAGES = {1054--1061},
      ISSN = {0022-314X,1096-1658},
   MRCLASS = {11D61 (11B68 11B83)},
  MRNUMBER = {2772487},
MRREVIEWER = {Pieter\ Moree},
       DOI = {10.1016/j.jnt.2011.01.004},
       URL = {https://doi-org.proxy2.library.illinois.edu/10.1016/j.jnt.2011.01.004},
}

@article {Moree2013,
    AUTHOR = {Moree, P.},
     TITLE = {Moser's mathemagical work on the equation
              {$1^k+2^k+\cdots+(m-1)^k=m^k$}},
   JOURNAL = {Rocky Mountain J. Math.},
  FJOURNAL = {The Rocky Mountain Journal of Mathematics},
    VOLUME = {43},
      YEAR = {2013},
    NUMBER = {5},
     PAGES = {1707--1737},
      ISSN = {0035-7596,1945-3795},
   MRCLASS = {11D61 (11A07)},
  MRNUMBER = {3127844},
MRREVIEWER = {Volker\ Ziegler},
       DOI = {10.1216/RMJ-2013-43-5-1707},
       URL = {https://doi-org.proxy2.library.illinois.edu/10.1216/RMJ-2013-43-5-1707},
}

@article {Moser1953,
    AUTHOR = {Moser, L.},
     TITLE = {On the diophantine equation {$1^n+2^n+3^n+\cdots +(m-1)^n=m^n.$}},
   JOURNAL = {Scripta Math.},
  FJOURNAL = {Scripta Mathematica},
    VOLUME = {19},
      YEAR = {1953},
     PAGES = {84--88},
      ISSN = {0036-9713},
   MRCLASS = {10.0X},
  MRNUMBER = {54627},
MRREVIEWER = {Ivan\ Niven},
}

@article {SondowMacMillan2011,
    AUTHOR = {Sondow, J. and MacMillan, K.},
     TITLE = {Reducing the {E}rd\H os-{M}oser equation
              {$1^n+2^n+\dots+k^n=(k+1)^n$} modulo {$k$} and {$k^2$}},
   JOURNAL = {Integers},
  FJOURNAL = {Integers. Electronic Journal of Combinatorial Number Theory},
    VOLUME = {11},
      YEAR = {2011},
     PAGES = {A34, 8},
      ISSN = {1553-1732},
   MRCLASS = {11D61},
  MRNUMBER = {2798610},
MRREVIEWER = {Attila\ B\'erczes},
       DOI = {10.1515/INTEG.2011.058},
       URL = {https://doi-org.proxy2.library.illinois.edu/10.1515/INTEG.2011.058},
}

@article {SondowMacMillan2017,
    AUTHOR = {Sondow, J. and MacMillan, K.},
     TITLE = {Primary pseudoperfect numbers, arithmetic progressions, and
              the {E}rd\H os-{M}oser equation},
   JOURNAL = {Amer. Math. Monthly},
  FJOURNAL = {American Mathematical Monthly},
    VOLUME = {124},
      YEAR = {2017},
    NUMBER = {3},
     PAGES = {232--240},
      ISSN = {0002-9890,1930-0972},
   MRCLASS = {11D68 (11A41)},
  MRNUMBER = {3626244},
MRREVIEWER = {Johnny\ Edwards},
       DOI = {10.4169/amer.math.monthly.124.3.232},
       URL = {https://doi-org.proxy2.library.illinois.edu/10.4169/amer.math.monthly.124.3.232},
}

@article {GrauOllerMarcen2022,
    AUTHOR = {Grau, J. M. and Oller-Marc\'en, A. M. and Sadornil, D.},
     TITLE = {On {$\mu$}-{S}ondow numbers},
   JOURNAL = {Acta Math. Hungar.},
  FJOURNAL = {Acta Mathematica Hungarica},
    VOLUME = {168},
      YEAR = {2022},
    NUMBER = {1},
     PAGES = {217--227},
      ISSN = {0236-5294,1588-2632},
   MRCLASS = {11D68 (11A07 11A51)},
  MRNUMBER = {4517540},
MRREVIEWER = {Joshua\ Zelinsky},
       DOI = {10.1007/s10474-022-01271-w},
       URL = {https://doi-org.proxy2.library.illinois.edu/10.1007/s10474-022-01271-w},
}

@article {Erdos1949,
    AUTHOR = {Erd\H{o}s, P.},
     TITLE = {Advanced Problem 4347},
   JOURNAL = {Amer. Math. Monthly},
  FJOURNAL = {American Math Monly},
    VOLUME = {56},
      YEAR = {1949},
     PAGES = {343},
      ISSN = {},
   MRCLASS = {},
  MRNUMBER = {},
MRREVIEWER = {},
}

@article {BJM2000,
    AUTHOR = {Butske, W. and Jaje, L. M. and Mayernik, D. R.},
     TITLE = {On the equation {$\sum_{p|N}(1/p)+(1/N)=1$}, pseudoperfect
              numbers, and perfectly weighted graphs},
   JOURNAL = {Math. Comp.},
  FJOURNAL = {Mathematics of Computation},
    VOLUME = {69},
      YEAR = {2000},
    NUMBER = {229},
     PAGES = {407--420},
      ISSN = {0025-5718,1088-6842},
   MRCLASS = {11D68 (05C50 11Y50)},
  MRNUMBER = {1648363},
MRREVIEWER = {Pieter\ Moree},
       DOI = {10.1090/S0025-5718-99-01088-1},
       URL = {https://doi-org.proxy2.library.illinois.edu/10.1090/S0025-5718-99-01088-1},
}

@article {Moree2011,
    AUTHOR = {Moree, P.},
     TITLE = {A top hat for {M}oser's four mathemagical rabbits},
   JOURNAL = {Amer. Math. Monthly},
  FJOURNAL = {American Mathematical Monthly},
    VOLUME = {118},
      YEAR = {2011},
    NUMBER = {4},
     PAGES = {364--370},
      ISSN = {0002-9890,1930-0972},
   MRCLASS = {11D41},
  MRNUMBER = {2800348},
       DOI = {10.4169/amer.math.monthly.118.04.364},
       URL = {https://doi-org.proxy2.library.illinois.edu/10.4169/amer.math.monthly.118.04.364},
}

\end{document}